\documentclass[11pt]{article}
\usepackage{amsmath,amsthm,amsfonts,latexsym,amssymb,enumerate,color}
\usepackage{graphicx}

\def\CC{\mathbb C}
\def\DD{\mathbb D}
\def\RR{\mathbb R}

\def\w{\omega}

\def\LL{\mathcal L}

\def\epsilon{\varepsilon}

\def\imag{\mathop{\rm Im}\nolimits}
\def\re{\mathop{\rm Re}\nolimits}
\def\Re{\re}

\def\dist{\mathop{\rm dist}\nolimits}

\def\beq{\begin{equation}}
\def\eeq{\end{equation}}

\def\beginpf{\begin{proof}}
\def\endpf{\end{proof}}

\newtheorem{thm}{Theorem}[section]

\newtheorem{lem}[thm]{Lemma}
\newtheorem{cor}[thm]{Corollary}
\newtheorem{rem}[thm]{Remark}

\begin{document}

\title{Laplace--Carleson embeddings on model spaces and boundedness of truncated Hankel and Toeplitz operators}

\author{Jonathan R. Partington\thanks{School of Mathematics, University of Leeds, Leeds LS2 9JT, UK.
\tt j.r.partington@leeds.ac.uk}, \ 
 Sandra Pott\thanks{Faculty of Science,
Centre for Mathematical Sciences,
Lund University,
22100 Lund, Sweden. \tt sandra.pott@math.lu.se},
\ and Rados\l aw Zawiski\thanks{Institute of Automatic Control, Silesian University of Technology, 44-100 Gliwice, Poland. 
\tt  	radoslaw.zawiski@polsl.pl}\ \thanks{The research presented here was done while the third author was a Marie Curie Research Fellow at the School of Mathematics, University of Leeds, UK.}
}

\maketitle

\begin{abstract}
A characterisation is given of bounded embeddings
from weighted $L^2$ spaces on bounded intervals into $L^2$ spaces on the half-plane, induced
by isomorphisms given by the Laplace transform onto  
weighted Hardy and Bergman spaces (Zen spaces). As an application necessary and
sufficient conditions are given for the
boundedness of truncated Hankel and Toeplitz integral operators, including the weighted case, and on model spaces.
\end{abstract}


\noindent {\bf 2010 Subject Classification:} 30H10, 32A36, 44A10, 47B35, 93B28

\section{Introduction}

By a Laplace--Carleson embedding, we mean the composition of two mappings: first, the Laplace
transform from a space $Z$ of functions defined on $(0,\infty)$ to a space $A$ of analytic functions
defined on the right-half plane $\CC_+$, then the natural embedding from $A$ into
$L^2(\CC_+,\mu)$ where $\mu$ is a positive regular Borel measure. As we shall make more
precise in Section \ref{sec:2}, the boundedness of such embeddings (which has applications in
the study of weighted Hankel operators \cite{aolo}, and linear systems \cite{JPP1,JPP2}) can be tested
by considering the measure of  Carleson squares. \\

The simplest example of the above situation is the case $Z=L^2(0,\infty)$, which, by the Paley--Wiener
theorem    is isomorphic to the Hardy space $H^2(\CC_+)$ by means of the Laplace transform.
The boundedness of the Laplace--Carleson embedding is then described by means of the
classical Carleson embedding theorem (see Section \ref{sec:2}).\\

Finite-interval embeddings, on the other hand, form a special case of an  embedding of a model space (an invariant subspace for
the backwards shift), for which a theory of Carleson measures in the unweighted case has been
studied by many authors, see e.g. \cite{cohn, baranov, LSSUW}.

Here we prove a very general theorem, Theorem \ref{thm:maince}, for a class of spaces (Zen spaces),
including the Hardy and Bergman spaces, for finite-interval embeddings. We also prove embedding theorems for a wider class of model spaces in case that $\mu$ is supported on a sector or on the real line (Thm. \ref{radialplane}, \ref{sectorplane}).\\

One application of these embedding results lies in the theory of boundedness of weighted
Hankel integral operators, including the Glover operator studied in rational $L^2$
approximation \cite{GLP}. This can even be used to characterize the
boundedness of certain truncated Hankel and Toeplitz operators acting on general model spaces (Thm. \ref{thmtrunc}, Rem. \ref{remtoep}).
We also discuss further applications in the theory of admissibility of control operators for linear 
semigroup systems 
(see, e.g. \cite{TW}).

\section{Background} \label{sec:2}

\subsection{Laplace--Carleson embeddings}

For $y \in \RR$ and $h>0$ the Carleson square 
$Q_{a,h}$ denotes the set 
\[
\{x+iy \in \CC_+: 0<x \le h, \ a-h/2 \le y \le a+h/2 \}.
\]
As usual we define the Laplace transform $\LL$, for suitable functions $f:(0,\infty) \to \CC$ (to be specified later)
via the formula
\[
(\LL f)(s)=\int_0^\infty f(t)e^{-st} \, dt, \qquad (s \in \CC_+).
\]
Recalling the Paley--Wiener theorem \cite[Chap. 19]{rudin}
that the Laplace transform provides an isomorphism of $L^2(0,\infty)$ onto
$H^2(\CC_+)$, the Hardy space of the right-hand half-plane $\CC_+$,
we may write the
 standard form of the Carleson embedding theorem (see, e.g., \cite[Thm. I.5.6]{garnett}) as follows.

\begin{thm}
Let $\mu$ be a positive regular Borel measure on $\CC_+$.
Then the following are equivalent.
\begin{enumerate}
\item The Laplace--Carleson embedding $\LL: L^2(0,\infty) \to L^2(\CC_+, \mu)$ is well-defined and bounded.
\item There is a constant $C_1>0$ such that for each $w \in \CC_+$
\[
\int_{\CC_+} \frac{d \mu(s)}{|s+w|^2} \le \frac{C_1}{\re w}.
\]
\item There is a constant $C_2>0$ such that $\mu(Q_{a,h}) \le C_2\,h$ for each Carleson
square in $\CC_+$.
\end{enumerate}
Moreover, the norm of the embedding and the constants $C_1$ and $C_2$ are all equivalent.
\end{thm}

The above theorem was generalized considerably in \cite{JPP1}, 
with applications to problems in control theory given in \cite{JPP2}.
The context used was
that of Zen spaces, defined using a measure
$\nu = \tilde\nu \otimes \lambda$, where $\tilde\nu$ is a positive regular measure on $[0,\infty)$ with $0 \in \mathrm{supp} \tilde \nu$,
satisfying the $\Delta_2$ condition 
\begin{equation}   \label{deltacond}
\sup_{t > 0} \frac{\tilde\nu[0,2t)}{\tilde\nu [0,t)} < \infty,
\end{equation}
and $\lambda$ denotes Lebesgue
measure. Defining
\beq\label{eq:defw}
w(t)=2\pi \int_0^\infty e^{-2rt} d\tilde\nu(r) \qquad (t>0)
\eeq
there is an isometric map $\LL: L^2(0,\infty; w(t)\, dt) \to A^2_\nu$,
where
\[
A^2_\nu= \left\{ f: \CC_+ \to \CC \hbox{ analytic }: \sup_{\varepsilon>0} \int_{\overline{\CC_+}} |f(z+\epsilon)|^2 \, d\nu(z) < \infty \right \}
\]
(see \cite[Prop. 2.3]{JPP1} and the references therein).
Note that $\tilde\nu= \delta_0$, a Dirac measure, corresponds to the Hardy space and $w(t)=2\pi$ for all $t$, 
while $\tilde\nu= \lambda$ corresponds
to the Bergman space with $w(t)=\pi/t$. 

The following result, which we shall extend to finite-time embeddings, is given in \cite[Thm. 2.1]{JPP1}.

\begin{thm}\label{thm:2.2r}
With $\nu$ and $w$ as above, the Laplace transform $\LL$ defines a bounded linear map
from $L^2(0,\infty; w(t)\, dt)$ into $L^2(\CC_+,\mu)$ 
if and only if there is a constant $k>0$ such that
for each Carleson square $Q_{a,h}$   we have
\beq\label{eq:creal21}
\mu(Q_{a,h}) \le
k \nu(Q_{a,h})  .
\eeq
\end{thm}

\subsection{Model spaces on the disc and half-plane}

For an inner function $\phi \in H^\infty(\DD)$, with $\DD$ the unit disc, we write $K^\DD_\phi= H^2(\DD) \ominus \phi H^2(\DD)$,
the associated model space. These are the nontrivial subspaces invariant under the backward
shift operator, and have been much studied. A key reference is \cite{GMR}.

Likewise, for an inner function $\theta \in H^\infty(\CC_+)$,  
we have $K^{\CC_+}_\theta =  H^2(\CC_+) \ominus \theta H^2(\CC_+)$. 
We may use the self-inverse bijection $M$ with $M(z)=(1-z)/(1+z)$ between $\DD$ and ${\CC_+}$,
and relate an inner function $\phi \in H^\infty(\DD)$ with an 
inner function $\theta \in H^\infty(\CC_+)$ by $\theta=\phi \circ M$.

We will be particularly interested in the  case that $\theta_T(s)=\exp(-Ts)$ for some $T>0$, and thus
$\phi_T(z)=\exp(T(z-1)/(z+1))$. This is because $L^2(0,T)$ is mapped isomorphically by the
Laplace transform onto the Paley--Wiener model space $K^{\CC_+}_{\theta_T}$.

In the context of Carleson embeddings from model spaces we may construct the following commutative
diagram:

\begin{equation}\label{eq:comdi}
\begin{array}{ccc}
K^\DD_\phi & \hookrightarrow & L^2(\DD,\rho) \\
\downarrow && \downarrow \\
K^{\CC_+}_\theta & \hookrightarrow & L^2(\CC_+,\mu)
\end{array}
\end{equation}
where  the relation between the measures $\rho$ and $\mu$ is
\beq\label{eq:defrho}
d\mu(s)=|1+s|^2\,d\rho(z)=4  \, d\rho(z)/|1+z|^2,
\eeq
where $s = M(z) = (1-z)/(1+z)$.
Indeed $F \in K^{\CC_+}_\theta$ if and only if $f \in K^\DD_\phi$, where
\[
f(z)=F(M(z))/(1+z),
\]
 since to within a constant the mapping from $F$ to $f$ is an isometry
from $H^2(\CC_+)$ onto $H^2(\DD)$. It takes $\theta H^2(\CC_+)$
onto $\phi H^2(\DD)$, and   preserves
orthogonal complements. 

Using the mappings above, one can translate the results of Cohn \cite{cohn}, Lacey et al \cite{LSSUW} and others
into the half-plane for a wide class of inner functions. This is discussed  in Section \ref{sec:3.2}.

\section{Embedding theorems}
\subsection{Carleson embeddings for Zen spaces}

We now consider finite-interval embeddings induced by the Laplace transform. It will be seen that
the results hold for all $T$ (with the implied norm bounds depending on $T$).
We write $\|.\|_w$ for the norm in $L^2(0,\infty;w(t) \, dt)$.

\begin{thm}\label{thm:maince}
With $\nu$ and $w$ related as in \eqref{eq:defw} and $T>0$, the Laplace transform $\LL$ defines a bounded linear map
from $L^2([0,T], w(t)\, dt)$ into $L^2(\CC_+,\mu)$ 
if and only if there is a constant $k>0$ such that
for each Carleson square $Q_{a,h}$ with $h \ge 1$ we have
\beq\label{eq:crealfinite}
\mu(Q_{a,h}) \le
k \nu(Q_{a,h})  .
\eeq
Moreover, there is an equivalence between the
norm of the embedding and the constant $k$ in the Carleson measure condition.
\end{thm}

\beginpf
We define a new measure on $\CC_+$ by $\mu_1(E)=\mu(E-1)$,
so that $\mu_1$ is a shifted version of $\mu$, supported on the half-plane $\{z \in \CC: \re z > 1 \}$.
For $f \in L^2(0,T; w(t)\, dt)$ we define $f_1(t)=e^t f(t)$, so that $f_1 \in L^2(0,T;w(t)\, dt)$
with $\|f_1\|_w \le e^T \|f\|_w$.
Note that $\LL f_1(s)=\LL f(s-1)$, which is well-defined as $\LL f$ is an entire function.

By (\ref{eq:crealfinite}), we have that $\mu_1$ satisfies the  Carleson condition (relative to $\nu$), and so,
as in Theorem~\ref{thm:2.2r},
\[
\|\LL f\|_{L^2(\mu)} = \|\LL f_1\|_{L^2(\mu_1)} \le c \|f_1\|_w \le ce^T \|f\|_w,
\]
for some absolute constant $c$, and so the Laplace--Carleson embedding is bounded.

For the converse, let $\mu_1$ be defined as above and suppose that the embedding $ L^2([0,T],w(t)\, dt) \rightarrow L^2(\CC_+,\mu)$ is bounded.   We note that for each $\epsilon>0$, $w(t) \ge \tilde\nu[0,\epsilon) e^{-2\epsilon t}$ for all $t$.
Now we decompose $(0,\infty)$ into intervals $[nT/2,(n+1)T/2)$, for $n=0,1,2,\ldots$,
and consider
the norm of the embeddings $J_n:L^2([nT/2,(n+1)T/2),w(t)\, dt) \to L^2(\mu_1)$. Clearly  the embedding $J_1: L^2([T/2,T),w(t)\, dt) \to L^2(\mu_1)$ is bounded, since $\LL: L^2([0,T],w(t)\, dt) \rightarrow L^2(\mu)$ is bounded.

If $f_n \in  L^2([nT/2,(n+1)T/2),w(t)\, dt)$, then $\tilde f_n$ defined by $\tilde f_n(t)=f_n(t+(n-1)T/2)$ lies in
$L^2([T/2,T),w(t)\, dt)$ and
\[
\|\tilde f_n\|_w \le \left( \frac{w(T/2)}{w((n+1)T/2} \right)^{1/2} \|f_n\|_w,
\]
using the monotonicity of the weight $w$. Thus
\[
\|J_n\| \le \|J_1\|\left( \frac{w(T/2)}{\tilde\nu[0,\epsilon) e^{-2\epsilon (n+1)T/2}} \right)^{1/2} e^{-(n-1)T/2},
\]
since $\LL f_n(s)=\exp(-(n-1)Ts/2) \LL \tilde f_n(s)$.

Choosing $0<\epsilon < 1$, we see that $\sum_{n=1}^\infty \|J_n\|< \infty$, and thus
the embedding from
$L^2(0,\infty; w(t) \, dt)$
into $L^2(\CC_+, \mu_1)$ is bounded. Thus, using
Theorem 2.1 of \cite{JPP1} we see that $\mu_1(Q_{a,h}) \le
k \nu(Q_{a,h}) $ for all Carleson squares, which implies by (\ref{deltacond})  the same result for $\mu$, at least for $h \ge 1$ .

\endpf

From now on, we shall not state explicitly when there is an equivalence between norms of linear operators 
and the constants appearing in the test conditions, as this is generally clear.

\subsection{Embeddings of Model spaces}
\label{sec:3.2}

Recall the result of Cohn \cite{cohn} (a version of his Theorem 3.2):
let $\rho$ be a measure on $\overline\DD$ and write
\[
h(\xi)= \int_{\overline \DD} \frac{1-|\xi|^2}{|1-\overline\xi z| ^2} \, d\rho(z).
\]
The spectrum of an inner function $\phi$ consists of the closure
of the union of the zero set of $\phi$ and the support of the singular measure 
associated with $\phi$, if they exist. Equivalently it can be defined as
\[
\sigma(\phi)= \left\{ \zeta \in \overline{\DD}: \liminf_{z \to \zeta} |\phi(z)|=0 \right\}
\]
(see, for example, \cite[p.71]{nik1}).

\begin{thm}\cite[Thm. 3.2]{cohn}
Suppose $\phi$ is an inner function such that $D_\epsilon=\{z \in \DD: |\phi(z)|< \epsilon\}$
is connected for some $\epsilon>0$, and suppose that $\rho$ is a measure that assigns zero mass to     the
spectrum of $\phi$. 
Then the following are equivalent:
\begin{enumerate}
\item  \label{cohnbdd}  $\rho$ is a Carleson measure for $K^\DD_\phi$; that is, the natural (identity) mapping
$K_\phi^\DD \to L^2(\DD,\rho)$ is bounded;
\item \label{cohncond} there is a constant $c>0$ such that
\begin{equation}  \label{cohncrit}
h(\xi) \le \frac{c}{1-|\phi(\xi)|}
\end{equation}
for all $\xi \in \DD$.
\end{enumerate}
\end{thm}

The question of the equivalence of \ref{cohnbdd}.\ and \ref{cohncond}.\ is known as the Cohn conjecture. This was shown to fail for general model spaces by Nazarov and Volberg in 2002 \cite{naz};
however, it does hold for special cases. The advantage of the Cohn criterion (\ref{cohncrit}) is that it is fairly easy to check, compared to the test function conditions for the associate Clark measure used in \cite{LSSUW}, \cite{seip}. We will therefore be interested in settings in which versions of the Cohn conjecture hold. 
The following result by Cohn \cite{cohn2} for the case of radial embeddings will be very useful:
\begin{thm}\cite[Thm. 3]{cohn2} \label{radial}
Suppose $\phi$ is an inner function on $\DD$, and suppose that $\rho$ is a measure supported on $[-1, 1]$. 
Then the following are equivalent:
\begin{enumerate}
\item    $\rho$ is a Carleson measure for $K^\DD_\phi$; that is, the natural mapping
\begin{equation} \label{bdd}
K_\phi^\DD \to L^2(\DD,\rho) \text{ is bounded};
\end{equation}
\item there is a constant $c>0$ such that
\begin{equation}   \label{cond}
h(r) \le \frac{c}{1-|\phi(r)|}
\end{equation}
for all $r \in [-1, 1]$.
\end{enumerate}
Moreover, in case that $1 \notin \sigma(\phi)$ (resp. $-1 \notin  \sigma(\phi)$), it suffices to check condition (\ref{cond}) only for $r \le 1- \dist(1, \sigma(\phi))/2$ 
(resp.  $r \ge -1+ \dist(1, \sigma(\phi))/2$).
\end{thm}
\proof
The equivalence of (\ref{bdd}) and (\ref{cond}) is just Theorem 3 in \cite{cohn2}. 
The consideration of the case  $1 \notin \sigma(\phi)$ etc. is not part of Cohn's original statement, but it can easily be extracted from his proof: As in (16), (17) of \cite{cohn2}, reformulate 
 (\ref{cond}) for $0 \le r <1$ as
\begin{equation} \label{newtest}
(1 - |\phi(r)|) \int_{-1}^r \frac{1-r }{(1-t)^2} d \rho(t) + \frac{1- |\phi(r)|}{1-r} \rho([r,1]) \le C.
\end{equation}
Suppose that $1 \notin \sigma(\phi)$ and let $R = 1- \dist(1, \sigma(\phi))/2$.  Then $\phi$ extends analytically through $1$, and 
$$
 \frac{1- |\phi(r)|}{1-r} \le \tilde C  \frac{1- |\phi(R)|}{1-R} \text{ for } R \le r <1,
 $$
where $\tilde C$ is a constant only depending on $\phi$.
Hence for $R \le r <1$,
\begin{eqnarray*}
&& (1 - |\phi(r)|) \int_{-1}^r \frac{1-r }{(1-t)^2} d \rho(t) + \frac{1- |\phi(r)|}{1-r} \rho([r,1])\\
& \le & \tilde C  \frac{1 - |\phi(R)|}{ 1-R}  \left(    \int_{-1}^r \frac{(1-r )^2}{(1-t)^2} d \rho(t) +  \rho([r,1]) \right)\\
& \le & \tilde C  \frac{1 - |\phi(R)|}{ 1-R}  \left(    \int_{-1}^R \frac{(1-R )^2}{(1-t)^2} d \rho(t) +  \rho([R,1]) \right),
\end{eqnarray*}
where the last expression is the testing condition (\ref{newtest}) for $r=R$. The proof for $-1 \notin  \sigma(\phi)$ and $-1 < r \le 0$ follows similarly.
\qed

In the context of the equivalence between disc and half-plane given in \eqref{eq:comdi} and \eqref{eq:defrho},
we have
\begin{equation} \label{kernels}
 \int_{\DD} \frac{1-|\xi|^2}{|1-\overline\xi z| ^2} \, d\rho(z) = \int_{\CC_+} \frac{  4\re w}{|\overline w + s|^2} \, d\mu(s),
\end{equation}
where $\xi=(1-w)/(1+w)$ and $z=(1-s)/(1+s)$. \\

It is therefore possible to derive results on embeddings of spaces $K^{\CC_+}_\theta$ into
$L^2(\CC_+,\mu)$ from those in the setting of the unit disc.

In the spirit of Theorem \ref{thm:maince}, let us consider the case of the singular inner function  $\phi(z)= \exp(T(z-1)/(z+1))$ for some fixed $T>0$, which does satisfy the connectivity condition above, and whose spectrum is simply $\{-1\}$.
Then, using Cohn's result we deduce the following.

\begin{thm}
Suppose that $\mu$ is a regular Borel measure on $\overline{\CC_+}$ (assigning
zero mass to $\infty$).
Then the  Laplace transform maps $L^2(0,T)$ boundedly into $L^2(\CC+,\mu)$
if and only if there is a constant $c>0$ such that
\[
\int_{\CC_+} \frac{  \re w }{| w + s|^2} \, d\mu(s) \le \frac{c}{1-\exp(-T \re w)} \qquad (w \in \CC_+).
\]
\end{thm}

Note that   ${1-\exp(-T \re w)}$ tends to $1$ as $\re w \to \infty$ and
is asymptotic to $T \re w$ as $\re w \to 0$. It is clear that the condition does not depend on $T$,
which is as one would expect.
We therefore have a simpler condition for a bounded Laplace--Carleson embedding.

\begin{cor}\label{cor:14}
The Laplace transform maps $L^2(0,T)$ boundedly into\\ $L^2(\CC_+,\mu)$
if and only if there is a constant $c>0$ such that
\beq
\label{eq:cor02}
\int_{\CC_+} \frac{  d\mu(s)}{| w + s|^2}    \le \frac{c}{\re w}+\frac{c}{(\re w)^2} \qquad (w \in \CC_+).
\eeq
\end{cor}

This is the same as saying that the norm of the reproducing kernel $s \mapsto 1/(s+\overline w)$ in
$L^2(\CC_+,\mu)$ is dominated by its norm in the mixed Hardy  space with
norm $\|f\|=\left(\|f\|^2_{H^2}+ \|f\|^2_{H^\infty} \right)^{1/2}$.

From this we may obtain the following result, which is a special case of Theorem \ref{thm:maince},
and was first proved using the methods of this section before the more general theorem was derived.
We omit the details of the proof.

\begin{thm}\label{thm:CMms}
The Laplace transform maps $L^2(0,T)$ boundedly into\\ $L^2(\CC_+,\mu)$
if and only if there is a constant $k>0$ such that
for each Carleson square $Q_{a,h}$ with $h \ge 1$ we have
\beq\label{eq:creal5}
\nu(Q_{a,h}) \le
kh .
\eeq
\end{thm}
This is a special case of the Carleson embedding theorem for model spaces by Treil and Volberg in \cite{treil}.
In the case of general model spaces, Theorem \ref{radial} translates via the correspondences (\ref{eq:comdi}) and (\ref{kernels}) into
\begin{thm} \label{radialplane}
Suppose $\theta$ is an inner function on $\CC_+$, and suppose that $\mu$ is a measure supported on $[0, \infty)$. 
Then the following are equivalent:
\begin{enumerate}
\item    $\mu$ is a Carleson measure for $K^{\CC_+}_\theta$; that is, the natural mapping
\begin{equation} \label{planebdd}
K^{\CC_+}_\theta   \to L^2(\CC_+, \mu ) \text{ is bounded};
\end{equation}
\item there is a constant $c>0$ such that
\begin{equation}   \label{planecond}
\int_{0}^\infty \frac{  w }{| w + s|^2} \, d\mu(s)  \le \frac{c}{1-|\theta(w)|}
\end{equation}
for all $w \in (0, \infty)$.
\end{enumerate}
Moreover, in case that $0 \notin \sigma(\theta)$ (resp. $\infty \notin  \sigma(\theta)$), it suffices to check condition (\ref{planecond}) only for $w \ge \min\{ \dist(0, \sigma(\theta))/2, 1 \}$ 
(resp.  for $w \le     2 \sup\{|z| : z \in \sigma(\theta)\}$).
\end{thm}
\qed


\subsubsection{Carleson Embedding Theorems on model spaces for sectorial measures}
In the theory of infinite-dimensional linear systems, measures supported on a sector of  the right half plane appear frequently in a natural way, arising from the spectrum of the underlying semigroup. Embedding theorems for measures supported on such a sector are therefore of interests, and such embeddings behave in many ways better than in the general case, see e.g. \cite{JPP2}. To our knowledge, it is unknown whether
the Cohn conjecture holds for general model spaces in case of measures supported on a sector. Thm 1.1. and Section 2.5 in \cite{seip},  previous to \cite{LSSUW},  provide a relatively accessible criterion in terms of Clark measures and 
test functions in this case. The counterexample of Nazarov and Volberg \cite{naz} concerns measures supported on the unit circle (in the case of the disc) respectively the imaginary axis (in the case of the right half-plane), thus it does not provide information on the case of sectorial measures.
Here, we prove a version of the Cohn conjecture for sectorial measures and model spaces with some extra condition.
For the proof, we require the following  refinement of Lemma 6.1 in \cite{baranov}.

\begin{lem} \label{bara} Let $\theta$ be an inner function on $\CC_+$, let $i \w \in i \RR$ such that $i \w \notin \sigma(\theta)$, and let
$\Gamma$ be the cone  $ \{ z \in \CC_+: |\imag z- \w| \le \Re z \}$ . Then there exists a constant $C = C(\theta, \w) $ such that
$$
      \log |\theta(z)| \le C \frac{1}{\Re z}\frac{1}{|\theta(i\w + \delta/4)|}  |\theta'(i\w + \delta/4)|  \text{ for all } z = iw + s \in \Gamma, \quad s \ge \delta.
$$
where $\delta = \dist(i\w, \sigma(\theta))$.
\end{lem}
\proof We follow mainly the lines of the proof of Lemma 6.1 in \cite{baranov}.  By Frostman's Theorem, it suffices to verify the lemma for Blaschke products, so let us assume that $\theta$ is a Blaschke product with zeroes $(z_n)$. Since $\delta = \dist(i\w, \sigma(\theta)) > 0$, $|\theta|$ is bounded below on the Carleson square
$Q = \{i \w' + s : | \w - \w'| < \delta/4, 0 < s <\delta/2 \}$, and $\theta$ is continuous on the boundary of this square. 

Now let $z  \in \CC_+$, $\Re z \ge \delta$, $z \in \Gamma$. If $\theta(z) =0$, there is nothing to prove, so let $\theta(z) \neq 0$.
Note that for each $z_n$, 
\begin{equation*}
\begin{split}
& |i\w + \delta/4 + \bar z_n|^2\\
 =& |\w - \imag z_n|^2 + |\delta/4+ \Re z_n|^2 \ge \frac{\delta^2}{48 (\Re z)^2}\left(  |\imag z - \imag z_n|^2 + |\Re z + \Re z_n|^2 \right). \\
\end{split}
\end{equation*}
Then
\begin{equation*}
\begin{split}
  \log |\theta(z)| &= \frac{1}{2} \log |\theta(z)|^2 = 
     \frac{1}{2} \sum_{n} \log(\left|  \frac{z-z_n}{z + \bar z_n} \right|^2)=   \frac{1}{2}\sum_{n} \log\left(1-  \frac{4 \Re z_n \Re z}{|z + \bar z_n|^2} \right)\\
      &\le 
     - 2\Re z \sum_{n}  \frac{ \Re z_n}{|z + \bar z_n|^2} \\
     &\le 
     - \frac{\delta^2}{24} \frac{1}{\Re z}  \sum_{n}  \frac{ \Re z_n}{|i\w + \delta/4 + \bar z_n|^2} \\
      &\le 
     - \frac{\delta}{48} \frac{\delta}{\Re z}  \frac{ |\theta'(i\w + \delta/4)|}{|\theta(i\w + \delta/4)|} .\\
     \end{split}
\end{equation*}
Here, we use that for any $w \in \CC_+$ with $\theta(w) \neq 0$,
\begin{equation*}
\begin{split}
\left| \frac{\theta'(w)}{\theta(w)} \right| & = \left|  \sum_n \frac{ 2 \re z_n}{(w + \bar z_n)^2} \frac{w+ \bar z_n}{w-  z_n} \right|  =  \left|\sum_n \frac{ 2 \re z_n}{(w + \bar z_n) (w- z_n)} \right|\\
& \le  
   2  \sum_n  \frac{\re z_n}{|w + \bar z_n|^2 } \frac{|\bar w +  z_n|}{|\bar w-  \bar z_n|} . \\
 \end{split}                                            
\end{equation*}
\qed 

\begin{thm} \label{sectorplane}
Suppose $\theta$ is an inner function on $\CC_+$, and suppose that $\mu$ is a measure supported on a sector $\Gamma$ based at $0$. In addition, suppose that $0 \notin \sigma(\theta)$, and that either 
$\infty \notin  \sigma(\theta)$, or that $|\theta|$ is bounded away from $1$ on the intersection of $\Gamma$ with some right half plane.
Then the following are equivalent:
\begin{enumerate}
\item  \label{bddsec}  $\mu$ is a Carleson measure for $K^{\CC_+}_\theta$; that is, the natural mapping
\begin{equation} 
K^{\CC_+}_\theta   \to L^2(\CC_+, \mu ) \text{ is bounded};
\end{equation}
\item there is a constant $c>0$ such that
\begin{equation}  \label{cohnsect}
\int_{\CC_+} \frac{  w }{| w + s|^2} \, d\mu(s)  \le \frac{c}{1-|\theta(w)|}
\end{equation}
for all $w \in (0, \infty)$.
\end{enumerate}
\end{thm}

\proof
Note that by Lemma \ref{bara}, 
\begin{equation}  \label{lowbound}
 \inf \{1-|\theta(\w)| : \w \in \Gamma, \delta \le \Re w \le N \} >0 
\end{equation}
for any $N >0$. That means, for $w $ in this set, condition (\ref{cohnsect}) of Theorem  \ref{sectorplane} will be just the ordinary boundedness condition for reproducing kernels, and the ordinary
Carleson condition for Carleson squares of sidelength between $\delta/2$ and $N$ follows. If both $0$ and $\infty$ are not in $\sigma(\theta)$, this implies for suitable $N >0$ the boundedness of the embedding
(\ref{bddsec}) by the sufficiency part of the embedding theorem of Treil and Volberg in \cite{treil}, which holds for any model space. If 
$|\theta|$ is bounded away from $1$ on the intersection of $\Gamma$ with some right half plane $\{ z \in \CC: \Re z \ge N\}$, then we obtain the ordinary
Carleson condition for Carleson squares of sidelength greater or equal to $\delta/2$, and again, the embedding theorem of Treil and Volberg  \cite{treil} implies the boundedness of the embedding. In this latter case, 
(\ref{lowbound}) is a simple consequence of Harnack's Theorem.

The necessity of the Cohn condition (\ref{cohnsect}) holds as always, since it follows from boundedness of the embedding on the reproducing kernels of the model space. This finishes the proof.
\qed

\section{Applications}

\subsection{Truncated Hankel and Toeplitz operators}

In   \cite{aolo}, Carleson-measure techniques were developed to obtain information on the boundedness of
certain integral operators, including the Glover operator studied in rational $L^2$ approximation
\cite{GLP}.
Here we take this further, and consider truncated Hankel operators acting on $L^2(0,T)$, with symbol
$h$   supported on $(0,2T)$, given by
\beq\label{eq:tho}
\Gamma_h u(t) = \int_0^T h(t+ \tau) u(\tau) \, d\tau.
\eeq
Note that by using the unitary reversal operator $Rf(x)=f(T-x)$ on $L^2(0,T)$,
and defining $g(x)=h(T-x)$ so that $g$ is supported on $(-T,T)$,
we obtain
\beq\label{eq:tto}
R \Gamma_h u(t)= \int_0^T g(t-\tau) u(\tau) \, d\tau,
\eeq
which is a truncated Toeplitz operator acting on $L^2(0,T)$. 
Taking Laplace transforms and writing
$G=\LL g$ and $\theta_T(s)=e^{-sT}$, as before,
we see also that (if bounded) $R\Gamma_h$ is unitarily equivalent to
a ``standard" truncated Toeplitz operator 
\[
A^{\theta_T}_G: f \mapsto P_{K_{\theta_T}} Gf,
\]
 acting on the model space
$K_{\theta_T} \subset H^2(\CC_+)$.
Clearly questions of
boundedness of \eqref{eq:tho} and \eqref{eq:tto} are equivalent.
Note that some  related work on Toeplitz and Hankel operators acting on Paley--Wiener spaces
may be found in \cite{rochberg,smith}. Here we consider more general situations.

We recall Widom's theorem \cite{widom} (see also \cite[Chap. 2]{power}) which applies to the case that
$h$ is given as the Laplace transform of a positive Borel measure $\mu$
defined on $\RR_+$; such $h$ are sometimes called {\em diffusive}, and correspond to
linear systems that can be modelled using the heat equation \cite{montseny1}.
Widom's theorem asserts that
in this situation $\Gamma_h$ (with $T=\infty$) is bounded if and only
if there is a constant $C>0$ such that
  $\mu(0,x) \le Cx$ for all $x>0$.

We can now  state a finite-interval version of Widom's theorem.

\begin{thm}\label{thm:4.1}
Let $h$ be given as the Laplace transform of a positive Borel measure $\mu$
defined on $\RR_+$. Then the operator $\Gamma_h$ on $L^2(0,T)$ is
bounded if and only if the conditions of Theorem \ref{thm:CMms} hold.
That is,   $\mu(0,x) \le Cx$ for all $x \ge 1$.
\end{thm}

\beginpf
The proof is a  modification of an argument in \cite{aolo}.
Defining $Z_\mu: L^2(0,\infty) \to L^2(\RR_+,\mu)$ by
\[
Z_\mu f(s) = \int_0^T e^{-st} f(t) \, dt,
\]
we have, if $Z_\mu$ is bounded,
\begin{equation} \label{factor}
\langle Z_\mu f, Z_\mu g \rangle = \int_{\RR_+} \int_0^T e^{-st} f(t) \, dt \int_0^T e^{-s\tau} \overline{g(\tau)} \, d\tau \, d\mu(s) = \langle \Gamma_h f, g\rangle,
\end{equation}   
for $f,g \in L^2(0,T)$, since by definition
\[
h(t+\tau)= \int_{\RR_+} e^{-s(t+\tau)} \, d\mu(s) \quad \hbox{for} \quad t,\tau >0,
\]
and so $\Gamma_h$ is bounded. The converse follows on taking $f=g$.
\endpf

For example, taking $d\mu(x)=dx/\sqrt x$ we see that
the function $h(t)=1/\sqrt{t}$ defines a bounded Hankel operator on the space
$L^2(0,T)$, although not on $L^2(0,\infty)$, since the Carleson measure
condition is not satisfied uniformly for all $x>0$.

\begin{rem}
{\rm
In the case that $h$ is the Laplace transform (sometimes called Fourier--Borel transform) of a complex measure $\mu$ defined on $\CC_+$
the condition $ \mu(Q_{a,h}) \le
kh$ for $h \ge 1$ will be sufficient for boundedness, although not necessary in general.
}
\end{rem}

Some extensions   to weighted Hankel operators were given in \cite{aolo}, and we   give a
brief account of a
finite-time version of these results. Supposing once more that $h$ is the Laplace transform of a positive
Borel
measure $\mu$ on $\RR_+$, consider the integral operator $\Gamma_{h,w}$ defined on $L^2(0,T)$ by
\[
\Gamma_{h,w}u(t)=\int_0^T w(t) h(t+\tau) w(\tau) u(\tau) \, d\tau.
\]
As above, the boundedness of this operator is easily seen to be equivalent to the
boundedness of the Laplace transform from $L^2(0,T; dt/(w(t))^2)$ into $L^2(\RR_+,\mu)$,
and this is   a variation on Theorem \ref{thm:maince}.

As a special case, we consider the power weights, as in \cite[Thm.~2.7]{aolo}.

\begin{thm}\label{thm:boundedcm}
Let $w(t)=t^\alpha$ for $\alpha \ge 0$, and let $\mu \ge 0$ be a measure supported on $\RR_+$.
Then $\Gamma_{h,w}$ is bounded if and only if 
there is a $\gamma>0$ such that
i$\mu(0,x) \le \gamma x^{1+2\alpha}$ for all $x \ge 1$.
\end{thm}

\beginpf
The proof is similar to that of \cite[Thm.~2.7]{aolo}, the only change being that
we work with the interval $(0,T)$ rather than $(0,\infty)$;
boundedness of $\Gamma_{h,w}$ is equivalent to the boundedness of the Laplace--Carleson
embedding on $L^2(0,T; t^{-2\alpha} \, dt)$.  For $\alpha>0$,
this is characterised by
Theorem~\ref{thm:maince}, using the 
measure $d\tilde\nu(r)=r^{1+2\alpha} \, dr$ and the
formula
\[
\int_0^\infty e^{-2rt}r^{1+2\alpha} \, dr = \frac{\Gamma(2\alpha)}{2^{2\alpha}t^{2\alpha}}.
\]
The case $\alpha=0$  is covered by Theorem \ref{thm:4.1}.
\endpf


In the case of truncated Hankel operators on general model spaces, we can  proceed in exactly the same way as in Theorem \ref{thm:4.1}, using the embedding result, Theorem \ref{radialplane}. Let
$\theta: \CC_+ \rightarrow \CC$ be an inner function and
$K:=\LL^{-1}K_\theta \subseteq L^2(0, \infty)$ be the inverse Laplace transform of the corresponding model space. 
Then we obtain:
\begin{thm}\label{thmtrunc}
Let $h$ be given as the Laplace transform of a positive Borel measure $\mu$
defined on $\RR_+$ and let $\theta$ be an inner function on $\CC_+$. Then the truncated Hankel operator $\Gamma_h$ on $K$ given by
\begin{equation}\label{eq:thogen}
     \langle \Gamma_h f, g \rangle = \int_{\RR_+}  \int_{\RR_+} h(t + \tau) f(t) \overline{g(\tau)}  \, dt \, d \tau  \qquad (f,g \in K)
\end{equation}
is bounded, if and only if  condition (\ref{planecond}) of Theorem \ref{radial} holds, that is,
\begin{equation}   \label{hankelcond}
\int_{0}^\infty \frac{  w }{| w + s|^2} \, d\mu(s)  \le \frac{c}{1-|\theta(w)|} \quad \text{ for all }w \in (0, \infty).
\end{equation}
Moreover, in case that $0 \notin \sigma(\theta)$ (resp. $\infty \notin  \sigma(\theta)$), it suffices to check  (\ref{hankelcond}) only for $w \ge \min\{ \dist(0, \sigma(\theta))/2, 1 \}$ 
(resp.  for $w \le     2 \sup\{|z| : z \in \sigma(\theta)\}$).
\end{thm}
\qed

Furthermore,  a simple criterion for trace class  membership of $\Gamma_h$ is available in this setting. 
The formal adjoint of the 
Laplace--Carleson embedding
$$
   Z_\mu : K \rightarrow L^2(\CC_+, \mu), \quad Z_\mu f(s) = \int_0^\infty e^{-st} f(t) \, dt
$$
is given by 
$$
 Z_\mu^*: L^2(\CC_+, \mu) \rightarrow K, \quad Z_\mu^* g(t) = P_{K} \int_{\CC^+}  e^{-{\bar s}t} g(s) \, d\mu(s) =  \int_{\CC_+} u(s, t) g(s) \, d \mu(s),
$$
where the integral kernel $u(s,t)$ is given by
$$
   u(s,t) = P_{K} e^{- {\bar s} \cdot }(t) = k^\theta_s(t),
$$
and  $(k^\theta_s(t))_{s \in \CC_+}$ denotes the reproducing kernel of the model space $K$. Hence
\begin{multline}  \label{s2norm}
\| Z_\mu\|^2_{S_2} =  \int_{\CC_+} \| k^\theta_s\|^2 \, d\mu(s) 
=\frac{1}{4 \pi^2} \int_{\CC_+} \int_0^\infty \left|    \frac{1- {\bar \theta}(s) \theta( i \omega)}{\omega + \bar s} \right|^2 \, d \omega \, d\mu(s) \\
=\frac{1}{4 \pi^2} \int_{\CC_+}  \left|    \frac{1- |\theta(s)|^2}{2   s} \right|^2  d\mu(s).
\end{multline}
This identity was proved in \cite[Prop 5.1]{baranov} for the case of the unit disc.
Using the factorization identity (\ref{factor}) again, we obtain
\begin{thm}
Let $h$ be given as the Laplace transform of a positive Borel measure $\mu$
defined on $\RR_+$, let $\theta$ be an inner function on $\CC_+$ and $K_\theta $ the corresponding model space.  Then the trace norm of the truncated Hankel operator $\Gamma_h$ on $K=\LL^{-1}K_\theta$ is given by
$$
   \| \Gamma_h\|_{S_1} = \frac{1}{4 \pi^2} \int_{\CC_+}  \left|    \frac{1- |\theta(s)|^2}{2   s} \right|^2 \, d\mu(s).
$$
and $\Gamma_h$ is of trace class if and only if the right-hand side is finite.
\end{thm}

\begin{rem}  \label{remtoep}
As in the special case of \eqref{eq:tho} and \eqref{eq:tto}, the  truncated Hankel integral
operator $\Gamma_h$ as defined in \eqref{eq:thogen} is  equivalent to a 
truncated Toeplitz operator in the case when $\theta(\overline s)=\overline{\theta(s)}$ for $s \in i\RR$.

For it follows from \cite[Cor 4.4]{Pbook} that under these circumstances the Laplace transform of $\Gamma_h f$
equals $P_{K_\theta}( \LL h)( R\LL f)$, at least on a dense set of $f \in K$, where $R:L^2(i\RR) \to L^2(i\RR)$ denotes the reversal
operator $RF(s)=F(-s)$ ($s \in i\RR)$.

Now it easily verified that if $F \in  K_\theta$, then $\theta (R F)$ also lies in $ K_\theta$, and
has the same norm, at least for inner functions $\theta$ such that $\theta(\overline s)=\overline{\theta(s)}$.
Thus the operator $\Gamma_h$ in \eqref{eq:thogen} is  equivalent to a truncated
Toeplitz operator on $  K_\theta$ with symbol $\overline\theta (\LL h)$.
\end{rem} 

In general, results for other Schatten norms (even the Hilbert--Schmidt norm) are harder to
come by, as mentioned in \cite{CFT}. The interested reader may also consult \cite{LR}. 


\subsection{Admissibility of diagonal control systems}

Following the setting in \cite[Chap. 4, 5]{TW},
we take a $C_0$ semigroup $(T(t))_{t \ge 0}$  with
infinitesimal generator $A$, acting on a Hilbert space $X$, and consider the system
\[
\frac{dx(t)}{dt} = Ax(t)+Bu(t), \qquad x(0)=x_0, \qquad t \ge 0.
\]
Here $x(t) \in X$ is the state of the system at time $t$, $u(t) \in \CC$ the input and $B: \CC \to D(A^*)'$
the control operator.

Suppose now that $(T(t))_{t \ge 0}$ has
a Riesz basis of eigenvectors $(\phi_k)$, that is, $T(t)\phi_k = e^{\lambda_k t}\phi_k$ for each $k$,
where the sequence $(\lambda_k)$ lies in $\CC_-$.

We take a Hilbert space $Z$ of functions, which in \cite{JPP1,JPP2} is a weighted $L^2(0,\infty)$ space, but
in our case will be a weighted space $L^2(0,T;w(t) \, dt)$ for
some $T>0$. 

Let $B$ be a linear bounded map from $\CC$ to $D(A^*)'$, corresponding
to the sequence $(b_k)$. Then the control operator $B$ is said to be $Z$-admissible
for $(T(t))_{t \ge 0}$ if and only if there is a constant $m_0>0$ such that
\[ 
\left \| \int_0^\infty  T(t) Bu(t) \, dt\right\|_X \le m_0 \|u\|_Z, \qquad (u \in Z).
\]
As discussed in \cite{JPP1,JPP2}, $B$ will be
admissible if and only if the Laplace transform induces a bounded mapping from
$Z$ into $L^2(\CC_+,\mu)$, where $\mu$ is the measure $\sum |b_k|^2 \delta_{-\lambda_k}$.

Finally, the condition we arrive at is that it is
necessary and sufficient for $Z$-admissibility that there exist a constant $k>0$ such that
\[
\sum_{-\lambda_k \in Q_{a,h}} |b_k|^2 \le k \nu(Q_{a,h}) \qquad \hbox{for all Carleson squares with} \quad h \ge 1,
\]
where $\nu$ and $w$ are related is in (\ref{eq:defw}).
This is consistent with results produced by other methods, some of which are discussed in \cite{TW}.

\section{Acknowledgements}
We thank Roman Bessonov for a discussion concerning the reference \cite{cohn2}. The second author gratefully acknowledges support by Vetenskapsr{\aa}det, VR grant 2015-05552.
The work of the first and third authors has received funding from the European Union's Horizon 2020 research and innovation programme under the Marie Sk{\l}odowska-Curie grant agreement  no.\ 700833. Moreover, the research of the first author was supported by a guest professorship at Lund University, 
financed by the Wallenberg Foundation Program for Mathematics, KAW 2016.0436. \\

\end{document}